\documentclass[preprint,12pt]{elsarticle} 
\usepackage[english]{babel} 
 
\usepackage{enumerate}              
\usepackage{amsmath}   
\usepackage{amssymb}   
\usepackage{graphicx}      
\usepackage{epsfig}            
   
\usepackage{color}     
\usepackage{caption} 
\usepackage{float} 
\usepackage{tabulary} 

\usepackage{enumerate}      
\usepackage{amsmath}

\journal{...}  
\begin{document}

\makeatletter
\newcommand{\norm}[1]{\left\lVert#1\right\rVert}
\def\ln{\mathop{\operator@font ln}\nolimits}
\makeatother


\begin{abstract}
Based on the eigenvalue idea and the time-varying weighted vector norm in state space we construct here the lower and upper bounds on the solutions of uniformly asymptotically stable linear systems. We generalize the known results for the linear time-invariant systems to the linear time-varying ones.
\end{abstract} 

\title{Time-varying vector norm and lower and upper bounds on the solutions of uniformly asymptotically stable linear systems}
\author{R.~Vrabel}
\ead{robert.vrabel@stuba.sk}
\address{Slovak University of Technology in Bratislava, Institute of Applied Informatics, Automation and Mechatronics,  Bottova 25,  917 01 Trnava,   Slovakia}

\newtheorem{thm}{Theorem}
\newtheorem{lem}[thm]{Lemma}
\newtheorem{defi}[thm]{Definition}
\newdefinition{cor}{Corollary}
\newdefinition{rmk}{Remark}
\newdefinition{ex}{Example}
\newproof{pf}{Proof}
\newproof{pot1}{Proof of Theorem \ref{main}}
\newproof{pot2}{Proof of Theorem \ref{thm2}}

\pagestyle{headings}

\maketitle

\section[Introduction]{Introduction}

In addition to the Lyapunov stability criteria for the linear system of differential equations $\dot x=A(t)x,$ 
$\dot x=dx/dt,$ $t\geq t_0,$ $x\in\mathbb{R}^n,$ other types of conditions guaranteeing the stability often are useful. Typically these are sufficient conditions that are proved by application of the Lyapunov stability theorems \cite{Khalil}, or the Gronwall-Bellman inequality \cite{Chicone}, though sometimes either technique can be used, and sometimes both are used in the same proof of stability criterion. One of these theorems, providing the conditions for eventual stability of the linear systems is the following theorem.
\begin{thm}[\cite{Rugh}]\label{thm_rugh}
For the linear system $\dot x=A(t)x,$ $t\geq t_0$ denote the largest and smallest
point-wise eigenvalues of $A^T(t)+A(t)$ by $\lambda_{\max}(t)$ and $\lambda_{\min}(t).$ Then for any $t_0$ and $x(t_0)$
the solution $x(t)$ satisfies
\begin{equation}\label{thm_rugh_ineq}
\norm{x(t_0)}_I{\mathrm{e}}^{1/2\int\limits_{t_0}^t\lambda_{\min}(\tau)d\tau}\leq \norm{x(t)}_I\leq\norm{x(t_0)}_I{\mathrm{e}}^{1/2\int\limits_{t_0}^t\lambda_{\max}(\tau)d\tau},\ \ t\geq t_0.
\end{equation}
\end{thm}
Throughout the whole paper it is assumed that a matrix function $A(t):$ $[t_0,\infty)\to\mathbb{R}^{n\times n}$ is continuous.

This theorem belongs to the wider family of sufficient condition for stability of the linear systems based on the "logarithmic measure" of the system matrices \cite[p.~58, Theorem~3]{Coppel}.

Our aim in this paper is to prove more useful theorem based on the eigenvalues idea for estimating asymptotics of the solutions of uniformly asymptotically stable linear systems. The theory is illustrated by two examples. 

\subsection{Notations, definitions and preliminary results}
Let $\mathbb{R}^n$ denotes $n-$dimensional vector space over the real numbers, $x=(x_1,\dots,x_n)^T\in\mathbb{R}^n$ is a column vector and the  symbol $\norm{\cdot}$ refers to any (real) vector norm on $\mathbb{R}^n.$ Specifically, for a symmetric, positive definite real matrix $H,$ we define the weight $H$ vector norm $\norm{x}_H\triangleq \big( x^THx\big)^{1/2}.$ Obviously, for $H=I$ ($I=$ identity on $\mathbb{R}^n$) we obtain the Euclidean norm, $\norm{x}_I.$ For the matrices $H\in\mathbb{R}^{n\times n}$ as an operator norm we will use an induced norm. Particularly, for weight $H$ vector norm in $\mathbb{R}^n,$ the norm $\norm{M}_H=\big(\lambda_{\max}\big[\hat M^T\hat M\big]  \big)^{1/2}$ where $\hat M=H ^{1/2}MH^{-1/2},$ as was proved in \cite{Hu_Liu}. Further,  $\lambda_i\big[M\big],$ $i=1,\dots,n$ denotes the eigenvalues of the matrix $M$ and  $\lambda_{\min}\big[M\big]=\min\{\lambda_i\big[M\big]:\ i=1,\dots,n\}.$ 

In this paper we will deal solely with the uniformly asymptotically ($\Leftrightarrow$ uniformly exponentially) stable linear systems \cite[Theorem~4.11]{Khalil}, \cite[Theorem~6.13]{Rugh}; for the different types of stability and their relation, see e.~g. \cite{Zhou}. We say, that
\begin{defi}[\cite{, Khalil, Rugh}]
The linear system $\dot x=A(t)x$ is {\it uniformly asymptotically stable} (UAS) if
there exist finite positive constants $\gamma,$ $\lambda$ such that for any $t_0$ and $x(t_0)$ the corresponding solution satisfies
\[
\norm{x(t)}\leq\gamma\norm{x(t_0)}{\mathrm{e}}^{-\lambda (t-t_0)}, \ \ t\geq t_0.
\]
\end{defi}
\begin{thm}[\cite{, Khalil, Rugh}]\label{theorem_UAS}
The linear system $\dot x=A(t)x$ is uniformly asymptotically stable if and only if there exist finite positive constants $\gamma,$ $\lambda$ such that
\[
\norm{\Phi(t,\tau)}\leq\gamma{\mathrm{e}}^{-\lambda(t-\tau)}
\]
for all $t,\tau$ such that $t\geq\tau\geq t_0.$ The transition matrix $\Phi(t,\tau)\triangleq X(t)X^{-1}(\tau),$  where $X(t),$ $t\geq t_0$ is a fundamental matrix of the system $\dot x=A(t)x.$  If $A(t)=A,$ an $n \times n$ constant matrix, then the transition matrix $\Phi(t,\tau)={\mathrm{e}}^{A(t-\tau)}.$
\end{thm}
Theorem~\ref{thm_rugh} leads to proof of some simple criterion based on the eigenvalues of $A^T(t)+A(t);$ for a wider context in connection with so called "logarithm measure" of the matrices see also e.~g.  \cite{Afanasiev}, \cite{Dekker_Verwer}, \cite{Desoer}. 
\begin{cor}[\cite{Desoer} with \cite{Rugh}]
The linear system $\dot x=A(t)x$ is UAS if there exist finite positive constants $\tilde\gamma,$ $\tilde\lambda$ such that such that the largest point-wise eigenvalue of $A^T(t)+A(t)$ satisfies
\[
\bigg(2\ln\norm{\Phi(t,\tau)}_I\leq\bigg)\int\limits_{\tau}^t\lambda_{\max}\big[A^T(s)+A(s) \big]ds\leq\tilde\gamma-\tilde\lambda(t-\tau)
\]
for all $t,\tau$ such that $t\geq\tau\geq t_0.$ Then Theorem~\ref{theorem_UAS} will hold with $\gamma={\mathrm{e}}^{\tilde\gamma/2}$ and $\lambda=\tilde\lambda/2.$
\end{cor}
This criterion is quite conservative in the sense that many UAS linear systems do not satisfy the above condition as we now see.
\begin{ex}\label{example_LTI1}
The system $\dot x=Ax,$ $t\geq0$ with
\[
A=\left(
\begin{array}{cc} 
0 & \sqrt{10}\\ 
-\sqrt{10} & -2 
\end{array}
\right)
\]
is UAS because $\lambda_{1,2}\big[A\big]=-1\pm3\,{\mathrm{i}}.$ Then a straightforward computation and Theorem~\ref{thm_rugh} shows that $\norm{x(0)}_I{\mathrm{e}}^{-2t}\leq\norm{x(t)}_I\leq\norm{x(0)}_I$ for all $t\geq0.$
\end{ex}
Despite such examples the eigenvalue idea is not to be completely rejected. In Theorem~\ref{main} below we prove for the UAS linear systems $\dot x=A(t)x$ the stronger result than the inequality in Theorem~\ref{thm_rugh}.
\section{Main results}
The main results of this paper are summarized in the following theorem generalizing \cite[Theorem~3.1]{Hu_Liu} to the linear time-varying systems. Recall that although its claims are mainly of theoretical relevance, providing the necessary conditions for exponential stability, within its framework without giving details and exact mathematical explanation the important results regarding convergent systems were derived in \cite{Lohmiller}; for the definitions and comparisons with the notion of incremental stability see also \cite{Ruffer}. Moreover, this theorem provides also the lower bound on the solutions generally classified as difficult to obtain.
\begin{thm}\label{main}
Let the linear system $\dot x=A(t)x,$ with a continuous matrix function $A(t):$  $[t_0,\infty]\to\mathbb{R}^{n\times n}$ is UAS. Then there exists a continuous, symmetric and positive definite matrix function $H(t):$  $[t_0,\infty]\to\mathbb{R}^{n\times n}$ such that every solution $x(t)$ satisfies
\[
\left(\frac{\lambda_{\min}[H(t)]}{\lambda_{\max}[H(t)]}\right)^{1/2}\norm{x(t_0)}_I{\mathrm{e}}^{-\frac12\int\limits_{t_0}^t \frac{d\tau}{\lambda_{\min}[H(\tau)]}}\leq\norm{x(t)}_I
\]
\begin{equation}\label{main_ineq} 
\leq\left(\frac{\lambda_{\max}[H(t)]}{\lambda_{\min}[H(t)]}\right)^{1/2}\norm{x(t_0)}_I{\mathrm{e}}^{-\frac12\int\limits_{t_0}^t \frac{d\tau}{\lambda_{\max}[H(\tau)]}}\ \mathrm{for\ all}\ t\geq t_0,
\end{equation}
where
\begin{equation*}
H(t)=\int\limits_t^\infty\Phi^T(\tau,t)\Phi(\tau,t)d\tau,\qquad t\geq t_0
\end{equation*}
for non-constant system matrix $A(t)$,
\begin{equation*}
H=\int\limits_{0}^\infty {\mathrm{e}}^{A^T\tau}{\mathrm{e}}^{A\tau}d\tau
\end{equation*}
for constant system matrix $A$ and 
\begin{equation*}
\lambda_{\min}\big[H(t)\big]\leq\lambda_{\max}\big[H(t)\big]\leq\frac{\gamma^2}{2\lambda}.
\end{equation*}
Moreover, if $A(t)$ is bounded, $\norm{A(t)}_I\leq L$ for all $t\geq t_0$, then
\begin{equation}\label{thm_ineq2}
\frac1{2L}\leq\lambda_{\min}\big[H(t)\big]\leq\lambda_{\max}\big[H(t)\big]\leq\frac{\gamma^2}{2\lambda}.
\end{equation}
The positive constants $\gamma,$ $\lambda$ and the transition matrix $\Phi(t,\tau)$ are defined in Theorem~\ref{theorem_UAS}. 
\end{thm}
\begin{pf}
We begin with  the analysis of the properties of the matrix function $H(t),$ $t\geq t_0.$ Observe that $H(t)$ is symmetric and positive definite because such is the integrand $\Phi^T(\tau,t)\Phi(\tau,t)$ \cite[Corollary~14.2.10]{Harville}. The use of 
\begin{itemize}
\item the {\it Rayleigh-Ritz ratio} \cite{HornJohnson}, 
\item the fact that $\norm{\Phi(\tau,t)}_I=\norm{\Phi^T(\tau,t)}_I$ because every matrix and its transpose have the same characteristic polynomial \cite[Lemma 21.1.2]{Harville}, 
\item  the fact that spectral radius of the matrix $\Phi^T(\tau,t)\Phi(\tau,t)$ is less or equal to any induced matrix norm $\norm{\Phi^T(\tau,t)\Phi(\tau,t)},$ and 
\item  Theorem~\ref{theorem_UAS} 
\end{itemize}
yields for every fixed $t\geq t_0$ and $x\in\mathbb{R}^n$ that
\[
x^TH(t)x\leq\lambda_{\max}\bigg[\int\limits_t^\infty\Phi^T(\tau,t)\Phi(\tau,t)d\tau\bigg]\norm{x}^2_I
\]
\[
\leq \norm{\int\limits_t^\infty\Phi^T(\tau,t)\Phi(\tau,t)d\tau}_I\norm{x}^2_I
\]
\[
\leq \norm{x}^2_I\int\limits_t^\infty\norm{\Phi(\tau,t)}^2_Id\tau\leq\norm{x}^2_I\int\limits_t^\infty \gamma^2{\mathrm{e}}^{-2\lambda(\tau-t)}d\tau= \frac{\gamma^2}{2\lambda}\norm{x}^2_I.
\]
As a consequence, $\lambda_{\max}\big[H(t)\big]\leq\frac{\gamma^2}{2\lambda}$ because there is equality $x^TH(t)x=\lambda_{\max}\big[H(t)\big]\norm{x}^2_I$ for $x$ equal to the eigenvector corresponding to $\lambda_{\max}\big[H(t)\big].$
To prove the left inequality in (\ref{thm_ineq2}) we will need the following 
\begin{lem} Let $\norm{A(t)}_I\leq L$ for all $t\geq t_0.$
Then the solution $x(t)$ of the $\dot x=A(t)x$ satisfies
\begin{equation}\label{ineq_L}
\norm{x(t_0)}_I{\mathrm{e}}^{-L(t-t_0)}\leq \norm{x(t)}_I\leq \norm{x(t_0)}_I{\mathrm{e}}^{L(t-t_0)}, \ \ t\geq t_0.
\end{equation}
\end{lem}
Observe that the right-hand side inequality is uninteresting for UAS systems, every estimate of $\norm{x(t)}_I$ would grow exponentially as $t\to\infty.$
\begin{pf}
The claim of the lemma follows immediately from the chain of inequality
\[
\lambda_{\max}\big[A^T(t)+A(t)\big]\leq\norm{A^T(t)+A(t)}_I\leq 2\norm{A(t)}_I\leq 2L,
\]
\[
\lambda_{\min}\big[A^T(t)+A(t)\big]\geq-\norm{A^T(t)+A(t)}_I\geq -2\norm{A(t)}_I\geq -2L,
\]
and (\ref{thm_rugh_ineq}).
\end{pf}
Now let $\phi(\tau)$ is a solution of $d\phi/d\tau=A(\tau)\phi$ starting at $(t,x),$ that is, $\phi(\tau)=\Phi(\tau,t)x.$ Then for all $x\in\mathbb{R}^n$ 
\[
x^TH(t)x =x^T\bigg(\int\limits_t^\infty\Phi^T(\tau,t)\Phi(\tau,t)d\tau\bigg)x=\int\limits_t^\infty \phi^T(\tau)\phi(\tau)d\tau 
\]
and, by (\ref{ineq_L}),
\[
\int\limits_t^\infty\norm{\phi(\tau)}^2_Id\tau\geq \norm{x}^2_I\int\limits_t^\infty{\mathrm{e}}^{-2L(\tau-t)}d\tau= \frac1{2L}\norm{x}^2_I.
\]

Arguing analogously as above, $\lambda_{\min}\big[H(t)\big]\geq\frac1{2L}$ and the inequality (\ref{thm_ineq2}) is proved.

Now we are ready to prove the remaining part of the theorem, namely the inequality (\ref{main_ineq}).  Suppose $x(t)$ is a solution of $\dot x= A(t)x$ corresponding to a given $t_0$ and nonzero $x(t_0).$ 
Let us formally consider a time-varying weighted vector norm of the solutions $\norm{x(t)}_{H(t)}.$ Then
\[
\frac{d}{dt}\norm{x(t)}^2_{H(t)}=\frac{d}{dt}\bigg[ x^T(t)H(t)x(t)\bigg]
\]
\begin{equation}\label{pf_1}
=x^T(t)\bigg[A^T(t)H(t)+\dot H(t) +H(t)A(t)\bigg]x(t).
\end{equation}
Now we show that the function $H(t)$ satisfies
\begin{equation*}
\dot H(t)+A^T(t)H(t)+H(t)A(t)=-I. 
\end{equation*}
Using that 
\[
\frac{d}{dt}\Phi(\tau,t)=-\Phi(\tau,t)A(t),
\]
\[
\frac{d}{dt}\Phi^T(\tau,t)=-A^T(t)\Phi^T(\tau,t)
\]
\cite[p.~70]{Coddington_Levinson}, \cite[p.~62]{Rugh}, respectively, and
\[
\Phi(\tau=\infty,t)=0 \ \ (\Leftarrow \mathrm{UAS}), \ \ \Phi(t,t)=I,
\]
we obtain that 
\[
\dot H(t)=\int\limits_t^\infty\Phi^T(\tau,t)\bigg[\frac{\partial}{\partial t}\Phi(\tau,t)\bigg]d\tau
\]
\[
+\int\limits_t^\infty\bigg[\frac{\partial}{\partial t}\Phi^T(\tau,t)\bigg]\Phi(\tau,t)d\tau - I
\]
\[
=-\int\limits_t^\infty\Phi^T(\tau,t)\Phi(\tau,t)d\tau\, A(t)
\]
\[
-A^T(t)\,\int\limits_t^\infty\Phi^T(\tau,t)\Phi(\tau,t)d\tau -I
\]
\[
=-A^T(t)H(t)-H(t)A(t)-I.
\]
Returning to (\ref{pf_1}), $\frac{d}{dt}\norm{x(t)}^2_{H(t)}=-\norm{x(t)}^2_I.$ Dividing through by $\norm{x(t)}_{H(t)}^2$ which is positive at each $t\geq t_0,$  the {\it Rayleigh-Ritz ratio} yields
\[
-\frac{1}{\lambda_{\min}\big[H(t)\big]}\leq\frac{\frac{d}{dt}\norm{x(t)}^2_{H(t)}}{\norm{x(t)}^2_{H(t)}}=-\frac{\norm{x}^2_I}{x^TH(t)x}\leq-\frac{1}{\lambda_{\max}\big[H(t)\big]}.
\]
Integrating from $t_0$ to any $t\geq t_0$ one gets
\[
-\int\limits_{t_0}^t \frac{d\tau}{\lambda_{\min}\big[H(\tau)\big]}\leq\ln\norm{x(t)}_{H(t)}^2-\ln\norm{x(t_0)}_{H(t)}^2\leq -\int\limits_{t_0}^t \frac{d\tau}{\lambda_{\max}\big[H(\tau)\big]}.
\]
Exponentiation followed by taking the nonnegative square root gives for all $t\geq t_0$ the inequality
\begin{equation}\label{main_ineq_not_conv}
\norm{x(t_0)}_{H(t)}{\mathrm{e}}^{-\frac12\int\limits_{t_0}^t\frac{d\tau}{\lambda_{\min}[H(\tau)]}}\leq\norm{x(t)}_{H(t)}\leq\norm{x(t_0)}_{H(t)}{\mathrm{e}}^{-\frac12\int\limits_{t_0}^t\frac{d\tau}{\lambda_{\max}[H(\tau)]}}.
\end{equation}

Finally using "norm conversion rule" between different weight $H_1$ and $H_2$ (recall $H_1,$ $H_2$ are symmetric and positive definite matrices)
\[
\frac{\lambda_{\min}[H_1]}{\lambda_{\max}[H_2]}\leq\frac{\norm{x}^2_{H_1}}{\norm{x}^2_{H_2}}=\frac{x^TH_1x}{x^TH_2x}\leq\frac{\lambda_{\max}[H_1]}{\lambda_{\min}[H_2]}\ \ \mathrm{for}\ \ x\neq0, 
\]
we obtain the inequality (\ref{main_ineq}).
\end{pf} 
\begin{rmk}
Combining \cite[Lemma~2.3, Theorem~2.1]{Hu_Liu} and \cite[p.~58, Theorem~3]{Coppel} we obtain
\[
\norm{x(t_0)}_{\widetilde H}{\mathrm{e}}^{-\int\limits_{t_0}^t\frac{d\tau}{\lambda_{\min}[\widetilde H]}}\leq\norm{x(t)}_{\widetilde H}\leq\norm{x(t_0)}_{\widetilde H}{\mathrm{e}}^{-\int\limits_{t_0}^t\frac{d\tau}{\lambda_{\max}[\widetilde H]}}
\]
which is a special case of (\ref{main_ineq_not_conv}) if $H(t)=\widetilde H/2.$  Observe that $\widetilde H$ in \cite{Hu_Liu} satisfies the Lyapunov equation $A^T\widetilde H+A\widetilde H=-2I.$ Thus, Theorem~\ref{main} represents generalization to the time-varying systems. Moreover, because $x(t)=\Phi(t,t_0)x(t_0),$ and from the properties of induced matrix norm we have
\[
\left(\frac{\lambda_{\min}[H(t)]}{\lambda_{\max}[H(t)]}\right)^{1/2}{\mathrm{e}}^{-\frac12\int\limits_{t_0}^t \frac{d\tau}{\lambda_{\min}[H(\tau)]}}\leq\norm{\Phi(t,t_0)}_I
\]
\[
\leq\left(\frac{\lambda_{\max}[H(t)]}{\lambda_{\min}[H(t)]}\right)^{1/2}{\mathrm{e}}^{-\frac12\int\limits_{t_0}^t \frac{d\tau}{\lambda_{\max}[H(\tau)]}}
\]
for $t\geq\tau\geq t_0.$ The general idea of the proof follows e.~g. the proof of \cite[Theorem~6.4, p.~100]{Rugh} and so the proof is omitted here. The last inequality generalizes \cite[Theorem~3.1]{Hu_Liu} to the linear time-varying systems. Moreover, we get also the lower bound on the solutions.
\end{rmk}
\section{Simulation results}
\begin{ex}[{\it Example~\ref{example_LTI1} revisited}]\label{example_LTI2}
Let us consider again the system from Example~\ref{example_LTI1}. The matrix exponential 
\[
{\mathrm{e}}^{At}=\frac{{\mathrm{e}}^{-t}}{3}\,
\left(
\begin{array}{cc}
 3\,\cos3\,t+\sin3\,t & \sqrt{10}\,\sin3\,t\\ 
-\sqrt{10}\,\sin3\,t & 3\,\cos3\,t-\sin3\,t 
\end{array}
\right)
\]
and the weight
\[
H=\int\limits_0^\infty {\mathrm{e}}^{A^T\tau} {\mathrm{e}}^{A\tau} d\tau
=\left(
\begin{array}{cc}
 3/5 & \sqrt{10}/20\\ 
\sqrt{10}/20 & 1/2
\end{array}
\right).
\]
The eigenvalues $\lambda_{\min}\big[H\big]=11/20 - \sqrt{11}/20,$ $\lambda_{\max}\big[H\big]=11/20 + \sqrt{11}/20$
and the inequality (\ref{main_ineq_not_conv}) becomes
\begin{equation}\label{example_LTI_ineq}
\norm{x(0)}_H{\mathrm{e}}^{-\frac{10t}{11-\sqrt{11}}}\leq\norm{x(t)}_H\leq\norm{x(0)}_H{\mathrm{e}}^{-\frac{10t}{11+\sqrt{11}}},
\end{equation}
where $\norm{x}_H=\big(3x_1^2/5 + \left(\sqrt{10}/10\right)x_1x_2 + x_2^2/2\big)^{1/2}.$ The result of simulation in the Matlab environment demonstrating effectiveness of the developed approach 
is depicted in Fig.~\ref{exampleLTI_solution}.

\centerline{}

\begin{figure}[H] 
\captionsetup{singlelinecheck=off}
   \centerline{\psfig{file=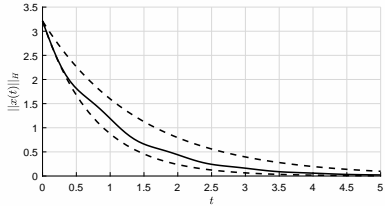,width=5.0cm} }
\caption{Solution of the linear time-invariant system from Example~\ref{example_LTI1} and~\ref{example_LTI2} with an initial state $x(0)=(x_1(0),x_2(0))^T=(-4,3)^T$ (the solid line) and the lower and upper bound given by 
(\ref{example_LTI_ineq}) (the dashed lines)}
\label{exampleLTI_solution}
\end{figure} 
\end{ex}
\begin{ex}\label{example_LTV}
For the linear time-varying system $\dot x=A(t)x,$ $t\geq0$ with 
\[
A(t)=\left(
\begin{array}{cc} 
-1 & {\mathrm{e}}^{-t}\\ 
0 & -3 
\end{array}
\right)
\]
the fundamental system (see, \cite{Zhou})
\[
X(t)=\left(\begin{array}{cc} {\mathrm{e}}^{-t} & \frac{{\mathrm{e}}^{-t}}{3}-\frac{{\mathrm{e}}^{-4\,t}}{3}\\ 0 & {\mathrm{e}}^{-3\,t} \end{array}\right).
\]
The eigenvalues of $A^T(t)A(t), t\geq0$
\[
\lambda_{1}\big[A^T(t)A(t)\big]=\frac{{\mathrm{e}}^{-2\,t}}{2}-\frac{{\mathrm{e}}^{-2\,t}}{2}\,\bigg[\left(4\,{\mathrm{e}}^{2\,t}+1\right)\,\left(16\,{\mathrm{e}}^{2\,t}+1\right)\bigg]^{1/2}+5 \to 1  
\]
as $t\to\infty,$ 
\[
\lambda_{2}\big[A^T(t)A(t)\big]=\frac{{\mathrm{e}}^{-2\,t}}{2}+\frac{{\mathrm{e}}^{-2\,t}}{2}\,\bigg[\left(4\,{\mathrm{e}}^{2\,t}+1\right)\,\left(16\,{\mathrm{e}}^{2\,t}+1\right)\bigg]^{1/2}+5 \to 9 
\]
as $t\to\infty;$ $\lambda_{1}\big[A^T(t)A(t)\big]<\lambda_{2}\big[A^T(t)A(t)\big]$ for all $t\geq0$ and 
$\norm{A(0)}_I=3.1796,$ $\norm{A(t)}_I=\big(\lambda_{\max}\big[A^T(t)A(t)\big]\big)^{1/2}\to 3$ (monotonically) as $t\to\infty$ and therefore the constant $L$ in (\ref{thm_ineq2}) is equal to $\norm{A(0)}_I=3.1796$ 
(Fig. \ref{exampleLTI_normAt}).

\centerline{}

\begin{figure}[H] 
\captionsetup{singlelinecheck=off}
   \centerline{\psfig{file=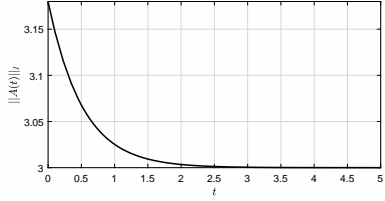,width=5.0cm} }
\caption{Time development of the $\norm{A(t)}_I,$ $t\geq0$}
\label{exampleLTI_normAt}
\end{figure} 
The transition matrix  
\[
\Phi(t,\tau)=X(t)X^{-1}(\tau)=\left(\begin{array}{cc} {\mathrm{e}}^{\tau-t} & \frac{{\mathrm{e}}^{-t}}{3}-\frac{{\mathrm{e}}^{3\,\tau-4\,t}}{3}\\ 0 & {\mathrm{e}}^{3\,\tau-3\,t} \end{array}\right)
\]
and the matrix function $H(t)$ from Theorem~\ref{main}
\[
H(t)=\int\limits_t^\infty \Phi^T(\tau,t)\Phi(\tau,t)d\tau=
\left(
\begin{array}{cc}
 \frac{1}{2} & \frac{{\mathrm{e}}^{-t}}{10}\\ 
\frac{{\mathrm{e}}^{-t}}{10} & \frac{{\mathrm{e}}^{-2\,t}}{40}+\frac{1}{6} 
\end{array}
\right)
\]
with the eigenvalues 
\begin{equation}\label{example_LTV_lmin}
\lambda_{\min}\big[H(t)\big]=\frac{{\mathrm{e}}^{-2\,t}}{80}-\frac{{\mathrm{e}}^{-2\,t}}{240}\,\bigg[336\,{\mathrm{e}}^{2\,t}+1600\,{\mathrm{e}}^{4\,t}+9\bigg]^{1/2}+\frac{1}{3}\to 1/6
\end{equation}
\begin{equation}\label{example_LTV_lmax}
\lambda_{\max}\big[H(t)\big]=\frac{{\mathrm{e}}^{-2\,t}}{80}+\frac{{\mathrm{e}}^{-2\,t}}{240}\,\bigg[336\,{\mathrm{e}}^{2\,t}+1600\,{\mathrm{e}}^{4\,t}+9\bigg]^{1/2}+\frac{1}{3}\to 1/2
\end{equation} as $t\to\infty$ (Fig.~\ref{exampleLTV_lambdas}). 

\centerline{}

\begin{figure}[H] 
\captionsetup{singlelinecheck=off}
   \centerline{
    \hbox{
     \psfig{file=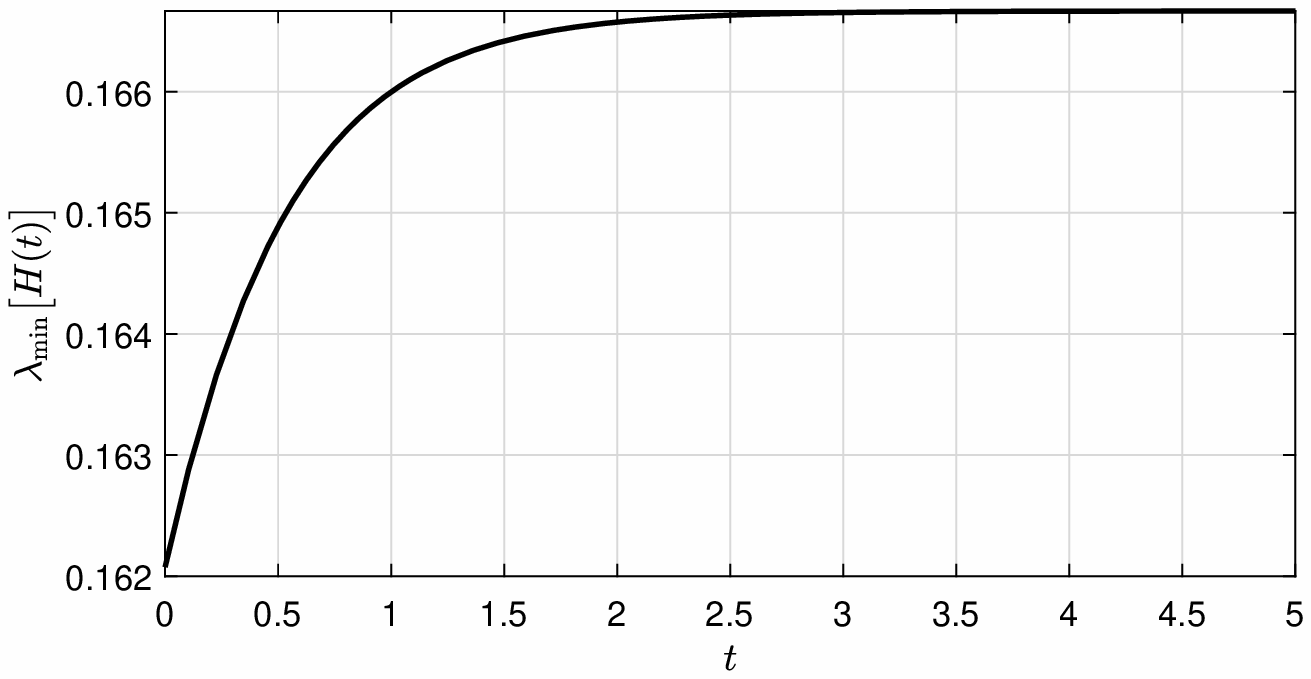,width=5.0cm, clip=}
     \hspace{1.cm}
     \psfig{file=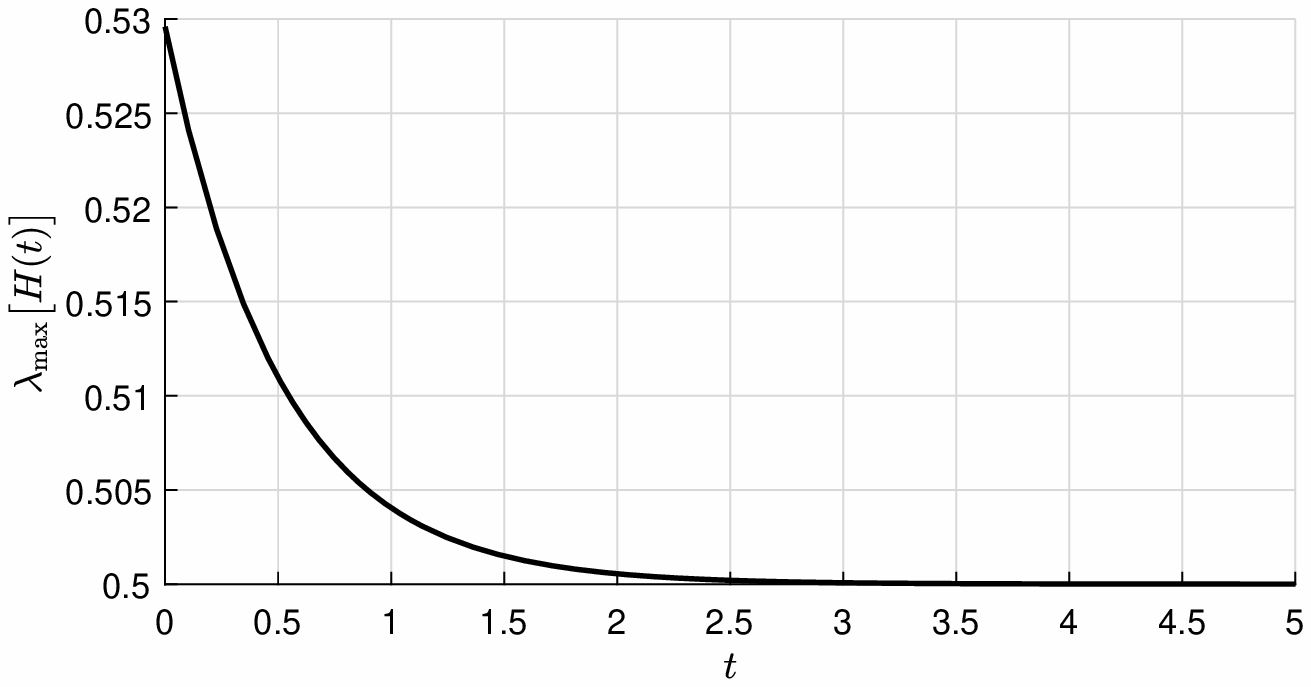,width=5.0cm,clip=}
    }
   }  
\caption{Time development of the functions $\lambda_{\min}(H(t))$ and $\lambda_{\max}(H(t)),$ $t\geq0$}
\label{exampleLTV_lambdas}
\end{figure} 
 The integrals in (\ref{main_ineq}) can be calculated explicitly
\[
-\frac{1}{2}\int\limits_0^t \frac{d\tau}{\lambda_{\min}\big[H(\tau)\big]}=\frac{3}{2}\,\ln\left(\rho-1\right)
-\frac{5}{2}\,\ln\left(\frac{2\,\sqrt{6}}{5}-\rho+\frac{7}{5}\right)
\]
\[
+\frac{1}{2}\ln\left(\left(\rho+1\right)\,\left(2\,\sqrt{6}-\rho+5\right)\right)+3.2375954052
\]
and
\[
-\frac{1}{2}\int\limits_0^t \frac{d\tau}{\lambda_{\max}\big[H(\tau)\big]}=\frac{3}{2}\,\ln\left(\rho+1\right)-\frac{5}{2}\,\ln\left(\frac{2\,\sqrt{6}}{5}+\rho+\frac{7}{5}\right)
\]
\[
+\frac{1}{2}\,\ln\left(\left(\rho-1\right)\,\left(2\,\sqrt{6}+\rho+5\right)\right)+2.1447615497,
\]
where 
\[
\rho = \left(\frac{100\,{\mathrm{e}}^{2\,t}+3\,\sqrt{6}+\frac{21}{2}}{100\,{\mathrm{e}}^{2\,t}-3\,\sqrt{6}+\frac{21}{2}}\right)^{1/2}.
\]
The result of simulation - the solution of system and lower and upper  bounds 
- are depicted in Fig.~\ref{exampleLTV_solution}.  

\centerline{}

\begin{figure}[H]  
\captionsetup{singlelinecheck=off}
   \centerline{\psfig{file=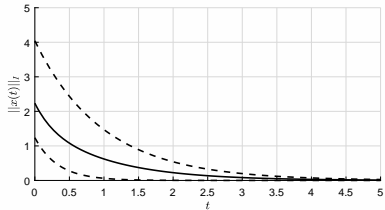,width=5.0cm} }
\caption{
Solution of linear time-varying system from Example~\ref{example_LTV} with an initial state $x(0)=(2,-1)^T$ (the solid line) and the lower and upper bound given by (\ref{main_ineq}), (\ref{example_LTV_lmin}) and 
(\ref{example_LTV_lmax}) (the dashed lines)}
\label{exampleLTV_solution}
\end{figure} 
Analyzing the properties of matrix function $H(t)$ it is obvious that
\[
\lambda_{\min}\big[H(0)\big]=\frac1{80}-\frac{\sqrt{1945}}{240}+\frac13(\approx 0.1621)\leq\lambda_{\min}\big[H(t)\big],
\]
\[
\lambda_{\max}\big[H(t)\big]\leq\lambda_{\max}\big[H(0)\big]=\frac1{80}+\frac{\sqrt{1945}}{240}+\frac13(\approx 0.5296)
\]
and
\[
(-t/2)\left(\lambda_{\min}\big[H(0)\big]\right)^{-1}=-3.0845t\leq-\frac{1}{2}\int\limits_0^t 1/\lambda_{\min}\big[H(\tau)\big]d\tau, 
\]
\[
(-t/2)\left(\lambda_{\max}\big[H(0)\big]\right)^{-1}=-0.9441t\geq -\frac{1}{2}\int\limits_0^t 1/\lambda_{\max}\big[H(\tau)\big]d\tau 
\]
for every $t\geq0.$ Thus we obtain more readable approximate estimate on the solutions
\[
0.5531\norm{x(0)}_I{\mathrm{e}}^{-3.0845t}\leq\norm{x(t)}_I\leq1.8075\norm{x(0)}_I{\mathrm{e}}^{-0.9441t}
\]
and Theorem~\ref{theorem_UAS} is satisfied for 
\[
\gamma=
\left(\frac{\lambda_{\max}[H(0)]}{\lambda_{\min}[H(0)]}\right)^{1/2}=\left(\frac{0.5296}{0.1621}\right)^{1/2}=1.8075,
\]
and 
\[
\lambda = (1/2)\left(\lambda_{\max}\big[H(0)\big]\right)^{-1}=0.9441.
\]

\end{ex}
\section*{Conclusion} 
In this paper we established the lower and upper bounds of {\bf all} solutions to uniformly asymptotically stable linear time-varying systems from the knowledge of {\bf one} fundamental matrix solution. Our approach is based on the eigenvalue idea and a time-varying metric on the state space $\mathbb{R}^n.$  The simulation experiments demonstrates the effectiveness of the proposed method for estimating solutions, generally classified as "difficult to obtain", especially in the case of the lower bounds.

\end{document}